\documentclass[11pt]{amsart}

\usepackage{amsmath, amsthm, amssymb}

\usepackage[T1]{fontenc}

\usepackage{comment}

\usepackage[abs]{overpic}		  % overpic automatically loads graphicx

\usepackage{xcolor}
\definecolor{indigo}{rgb}{0.29, 0.0, 0.51}  % custom colors
\usepackage[colorlinks, urlcolor=indigo, linkcolor=indigo, citecolor=indigo]{hyperref}

\theoremstyle{plain}
\newtheorem{theorem}{Theorem}

\newtheorem{question}[theorem]{Question}

\theoremstyle{definition}

\theoremstyle{remark}

\numberwithin{theorem}{section}
\def\co{\colon\thinspace}

\newcommand{\dfn}[1]{{\em #1}}        % definition
\newcommand{\R}{\mathbb{R}}           % the real numbers
\newcommand{\Z}{\mathbb{Z}}           % the integers
\makeatletter
\newcommand*\bigcdot{\mathpalette\bigcdot@{0.6}}
\newcommand*\bigcdot@[2]{\mathbin{\vcenter{\hbox{\scalebox{#2}{$\m@th#1\bullet$}}}}}
\makeatother

\DeclareMathOperator\tb{tb}                   % Thurston-Bennequin
\DeclareMathOperator\rot{rot}                 % rotation
\DeclareMathOperator\self{sl}                 % self linking
\begin{document}

% title
\title{Submanifolds in contact geometry} 

% author information
\author{John B. Etnyre}
\address{School of Mathematics \\ Georgia Institute of Technology \\ Atlanta \\ Georgia}
\email{etnyre@math.gatech.edu}

%\subjclass[2020]{57R17}

%abstract
\begin{abstract}
We survey what is known about various special types of submanifolds of contact manifolds and discuss their role in the development of contact geometry.
\end{abstract}

\maketitle
%\tableofcontents

%%%%%%%%%%%%%%%%%%%%%%%%%%%%%%%%%%%%
\section{Introduction}
%%%%%%%%%%%%%%%%%%%%%%%%%%%%%%%%%%%%

Studying submanifolds of contact manifolds has been a central part of the development of this vibrant subject. In this paper, we will survey what is known about such manifolds, the diverse tools one can use to study them, their impact on the general development of contact geometry, and current directions of research. In dimension $3$, where contact geometry is more fully developed, many of the milestones that have made the subject so rich arose out of the study and use of special submanifolds. For example:\\

\noindent
{\bf The existence of contact structures:} In \cite{Martinet71} and \cite{Lutz77}, Martinet and Lutz showed that any orientable $3$-manifold has a contact structure and, moreover, that any homotopy class of plane field on such a manifold was homotopic to a contact structure. Their construction involved surgeries on transverse knots in the standard contact structure $\xi_{std}$ on $S^3$. Since then, Ding and Geiges \cite{DingGeiges04} have shown that all contact structures on $3$-manifolds can be obtained by surgeries on Legendrian knots in $(S^3,\xi_{std})$, and Conway \cite{Conway19} has done the same for transverse knots. Analogous results were obtained in all dimensions by Conway and the author \cite{ConwayEtnyre20caps} and independently by Lazarev \cite{Lazarev2020caps}.

\noindent
{\bf Understanding types of contact structures:} In \cite{Bennequin83}, Bennequin showed that $\R^3$ had at least two distinct contact structures. One, the standard one, is characterized by an inequality on an invariant of transverse knots, while the other does not satisfy this inequality. This led Eliashberg to define the tight versus overtwisted dichotomy for contact structures on $3$-manifolds, \cite{Eliashberg89, Eliashberg92a}. These works of Bennequin and Eliashberg can be thought of as the birth of modern contact topology and have had far-reaching implications. %in low-dimensional topology.

\noindent
{\bf  Distinguishing tight contact structures:} Giroux \cite{Giroux1999} showed that the $3$-torus (and some other $3$-manifolds) has infinitely many distinct tight contact structures that were all homotopic as plane fields using the Legendrian knots contained within the tight contact structures. Since then, similar phenomena have been observed many times \cite{Giroux00, Honda00a}, and Van Horn-Morris and the author have shown that the ``transverse knot theory" a contact structure supports characterizes the contact structure \cite{EtnyreVanHorn-Morris11}. 

%While maybe not a milestone in the development of contact topology, we also note:

\noindent
%{\bf There is a rich theory of transverse and Legendrian knots:} 
{\bf Structure of Legendrian and transverse submanifolds:} Our understanding of transverse and Legendrian knots is still at an early stage, but we can already see many subtle and beautiful features. For example, there are many deep tools one can use to distinguish Legendrian knots that come from the theory of pseudo-holomorphic curves \cite{Chekanov02, Eliashberg1998, EliashbergGiventalHofer00}, from Heegaard Floer theory \cite{LiscaOzsvathStipsiczSzabo2009, OzsvathSzaboThurston2008}, and from micro-local sheaf theory \cite{CasalsZaslow2022, GuillermouKashiwaraSchapira2012, ShendeTreumannZaslow2017LegendrianSheaf}  (among others). The study of Legendrian submanifolds and invariants coming from pseudo-holomorphic curves in higher dimensions can be used to build a powerful invariant of smooth knots in $S^3$, \cite{EkholmEtnyreNgSullivan13}. We have classification theorems for simple knot types like the unknot \cite{EliashbergFraser09}, torus knots \cite{EtnyreHonda01b}, the figure eight knot \cite{EtnyreHonda01b}, cables of torus knots \cite{EtnyreHonda05, EtnyreLafountainTosun12}, and twist knots \cite{EtnyreNgVertesi13}. We also have various structure theorems about how the classification of Legendrian knots behaves under some cables \cite{ChatterjeeEtnyreMinMukherjee25nonloose}, some satellites \cite{EtnyreVertesi18}, and connected sum \cite{EtnyreHonda03}. All of these results are for transverse and Legendrian knots in tight contact manifolds, but there is also a surprisingly complex theory of such knots in overtwisted contact manifolds \cite{EliashbergFraser09, Etnyre13, EtnyreMinMukherjee22Pre, Matkovic22}. In forthcoming work \cite{EtnyreMinTosunVarvarezosPre}, the understanding of Legendrian knots in all contact structures on $S^3$ can also be used as a tool in classifying tight contact structures on certain manifolds. 

In the following sections, we will discuss each of the topics above in more detail and bring the reader up to the current state of research in each. We will end with a section discussing what is known about nice submanifolds of contact manifolds in higher dimensions. Due to space limitations, this survey will not be able to cover all aspects of the theory in detail, and so while we try to at least mention most aspects of the theory, the topics are skewed toward subjects the author knows best. 

We will forgo a review of basic definitions in contact geometry and refer the reader to any of the standard sources, such as \cite{Etnyre05, EtnyreTosunPre24, Geiges08}. 

\smallskip
\noindent
{\bf Acknowledgements.} I gratefully acknowledge the valuable input of Roger Casals, Lori Lejeune, and Bülent Tosun on early drafts of this paper.  The author was partially supported by National Science Foundation grant DMS-2203312 and the Georgia Institute of Technology’s Elaine M. Hubbard Distinguished Faculty Award.

%%%%%%%%%%%%%%%%%%%%%%%%%%%%%%%%%%%%%
%\section{Basic definitions and examples}
%%%%%%%%%%%%%%%%%%%%%%%%%%%%%%%%%%%%%
%
%Some info 
%
%include self-linking
%
%\begin{lemma}\label{tandlexistance}
%Any knot $K$ in a contact $3$-manifold $(M.\xi)$ may be $C^0$-approximated by a Legendrian knot and by a transverse knot. 
%\end{lemma}

%%%%%%%%%%%%%%%%%%%%%%%%%%%%%%%%%%%%
\section{Building contact structures}
%%%%%%%%%%%%%%%%%%%%%%%%%%%%%%%%%%%%
The first milestone in contact geometry we will discuss is the question of the existence of a contact structure on a given $3$-manifold $M$. More specifically, since a contact structure is a special plane field in the tangent space of $M$, we can ask whether every homotopy class of plane field $\eta$ admits a contact structure. Using some results from topology and a few basic facts about transverse knots in contact manifolds, one can remarkably answer this question.
\begin{theorem}[Martinet 1971, \cite{Martinet71} and Lutz 1977, \cite{Lutz77}]\label{originalexist}
If $\eta$ is a homotopy class of plane field on a compact, oriented manifold $M$, then there is a contact structure $\xi$ in the homtopy class $\eta$. 
\end{theorem}
We note that if a $3$-manifold admits a contact structure, then it must be oriented, so that hypothesis in the theorem is necessary. %We also remark that Martinet showed that any orientable $3$-manifold admited a contact structure, while Lutz provided the tools necessary to prove that there was a contact structure in every homotopy class of plane field (note that any $3$-manifold has infinitely many different homotopy classes of plane field). 
The main steps of the proof are as follows: 
(1) Any compact, oriented $3$-manifold can be obtained by Dehn surgery on some link $L$ in $S^3$ \cite{Lickorish62, Wallace60}. 
(2) The link $L$ is isotopic to a transverse link $T$. 
(3) The link $T$ has a standard neighborhood $N$. 
(4) Performing the Dehn surgery on $T$ removes this standard neighborhood $N$ and glues in solid tori. We have models of contact structures on solid tori that can be used to extend the contact structure on $\overline{S^3-N}$ to the solid tori, and thus to $M$. 
We now have one contact structure on $M$. Lutz introduced other models of contact structures on solid tori that could be used in Step~(4) above. This allows one to define a Lutz twist on a neighborhood of a transverse knot that does not change the ambient manifold but might change the contact structure. You can use Lutz twists to construct contact structures in all homotopy classes of plane fields. 

So we see that the geometric content (as opposed to the topological fact about Dehn surgery above) of our ``existence theorem" for contact structures relies on the existence of transverse realizations of knots, regular neighborhood theorems for transverse knots, and local models for surgery on transverse knots. These are fairly simple theorems, but they already have amazing consequences. 

Expanding on the results above, it turns out that all contact structures can be constructed by surgery on a Legendrian link. 
\begin{theorem}[Ding and Geiges 2004, \cite{DingGeiges04}]
Any contact structure on a compact $3$-manifold can be obtained from the standard tight contact structure on $S^3$ by contact surgery on a Legendrian link.
\end{theorem}
This theorem is a fundamental theorem that has many applications in the study of contact $3$-manifolds. We note that the theorem is closely related to the work of Honda and the author in \cite{EtnyreHonda02a}. %(in a similar way that Wallace's work is related to Lickorish's in the first item above). 
We also have an analogous theorem for surgery on transverse knots.
\begin{theorem}[Conway 2019, \cite{Conway19}]
Any contact structure on a compact $3$-manifold can be obtained from the standard tight contact structure on $S^3$ by admissible or inadmissible surgery on a transverse link.
\end{theorem}
The surgeries in this theorem are analogous to contact surgeries on Legendrian knots. The proofs of both these theorems rely on Eliashberg's beautiful classification of overtwisted contact structures discussed in the next section and the ideas in the proof of Theorem~\ref{originalexist}, but the main point of either of these theorems is that constructing any contact structure in dimension $3$ is closely related to understanding submanifolds of the standard contact structure on $S^3$.

%%%%%%%%%%%%%%%%%%%%%%%%%%%%%%%%%%%%
\section{Understanding properties of contact structures}
%%%%%%%%%%%%%%%%%%%%%%%%%%%%%%%%%%%%
After one knows that contact structures exist in a given homotopy class of plane field on an oriented $3$-manifold $M$, one is interested in whether there is more than one contact structure. The second milestone we will discuss is that there can be more than one, and they have very different properties. The story starts with work of Bennequin. In \cite{Bennequin83}, he showed that the contact structure 
\[
\xi_{std}= \ker (dz-y\, dx)
\]
and 
\[
\xi_{ot}=\ker (\cos r\, dz + r\sin r\, d\theta)
\]
on $\R^3$ (expressed in Euclidean and cylindrical coordinates, respectively) are distinct contact structures, up to contactomorphism. Specifically, he proved the following theorem.
\begin{theorem}[Bennequin 1983, \cite{Bennequin83}]
If $T$ is a transverse knot in $(\R^3,\xi_{std})$ that bounds a surface $\Sigma$, then \[\self(T)\leq -\chi(\Sigma),\] where $\chi(\Sigma)$ is the Euler characteristic of the surface. 
\end{theorem}
If one considers the disk in the $r\theta$-plane in $(\R^3,\xi_{ot})$ of radius a little larger than $\pi$, then it is easy to check that its boundary is a transverse knot with $\self=1$ and thus the contact structures $\xi_{std}$ and $\xi_{ot}$ must be different. So we have distinguished two contact structures based on the submanifolds one can find inside them!

The inequality in the theorem is now known as the Bennequin inequality. This work, and $h$-principle considerations, led Eliashberg to his revolutionary definition of tight and overtwisted contact structures \cite{Eliashberg89, Eliashberg92a}. A contact structure is \dfn{overtwisted} if it admits a Legendrian unknot with $\tb=0$, otherwise it is \dfn{tight}. (This is not Eliashberg's original formulation of overtwisted, but it is implied by it, and now we know it easily implies it.) Eliashberg also showed that in a tight contact structure, Bennequin's inequality always holds, and can be phrased in terms of Legendrian knots as follows: for a Legendrian knot $L$ that bounds a surface $\Sigma$, we have
\[
\tb(L)+|\rot(L)|\leq -\chi(\Sigma).
\]
Moreover, he showed that each homotopy class of plane field on an oriented $3$-manifold admits a unique overtwisted contact structure. This last result tells us that overtwisted contact structures are determined by algebraic topology; specifically, it is related to the obstruction theory associated with reducing the structure group of the tangent bundle from $SO(3)$ to the subgroup $SO(2)$. On the other hand, the Bennequin inequality indicates that tight contact structures can detect subtle information about the underlying $3$-manifold (such as bounding the minimal genus of a Seifert surface for a knot). 

The central place that the Bennequin inequality holds in contact geometry can be seen by the diverse ways in which it can be proven. The original, very topological proof, was through braid theory. Bennequin showed one can ``braid" any transverse knot in the standard contact structure, and then he analyzed ``braid foliations" on the knot's Seifert surface to prove his theorem. More recently, there have been proofs using pseudoholomorphic curves, convex surfaces, Seiberg-Witten theory, and Heegaard-Floer theory, to name just a few. 

We also note that due to work of Wu \cite{Wu06}, we can use Legendrian knots to detect overtwisted contact structures in a more subtle way. Specifically, recall that a small Seifert fibered space is an orbifold bundle over the $2$-sphere with three singular fibers. Wu showed that in some cases such a space must always have a regular fiber with a Legendrian representative along which the contact planes twist $0$ times (with respect to the product framing) and in other cases this implies the contact structure is overtwisted.  %and an important invariant of these spaces is denoted $e_0\in \Z$. The work in \cite{Wu06} shows that for a contact structure on a Seifert fibered space with $e_0\leq -3$, if the regular fiber of the Seifert fibered space admits a Legendrian representative along which the contact planes twist $0$ times (with respect to the product framing), then the contact structure is overtwisted. Interestingly enough, this is very sensitive to $e_0$. For example, if $e_0\geq 0$ then any tight contact structure will admit a Legendrian realization of a regular fiber with twisting $0$. 

%%%%%%%%%%%%%%%%%%%%%%%%%%%%%%%%%%%%
\section{Distinguishing tight contact structures}
%%%%%%%%%%%%%%%%%%%%%%%%%%%%%%%%%%%%
Transverse, or Legendrian, knots were first used to distinguish contact structures in the work of Bennequin, as discussed above, but there the knots distinguished between tight and overtwisted contact structures. The first time tight contact structures were distinguished by such knots was in work of Giroux. In \cite{Giroux1999}, he showed that any $T^2$-bundle over $S^1$, and in particular the $3$-torus $T^3$, admits an infinite number of tight contact structures. The result for $T^3$ was obtained earlier by Giroux and Kanda \cite{Kanda97}. To illustrate how one distinguishes the contact structures, we will focus on $T^3$ and think of it as the quotient $\R^3/\Z^3$ of $\R^3$ by the action of the integral lattice $\Z^3$. None of the contact structures 
\[
\xi_n=\ker(\cos 2\pi nz\, dx + \sin 2\pi nz\, dy),
\] 
on $T^3$, for $n$ a positive integer, are contactomorphic. One can show this by considering the knot $K$ given by the image of the $z$-axis in $\R^3$. When considering Legendrian representatives of $K$, one can show that the maximal possible twisting of the contact planes of $\xi_n$ along the knot is $-n$. 
%From this, it is easy to see that the $\xi_n$, for distinct $n$, are not contactomorphic. 
To prove the $\xi_n$ are not contactomorphic, one needs to consider Legendrian representatives of other knots. 
It is not hard to show that there are three linear knots in $T^3$ and the isotopy class (and contactomorphism class) of any tight contact structure on $T^3$ is determined by the Legendrian representatives of these knots. See \cite{Giroux00, Honda00a, Kanda97}. %(none of these papers precisely proved this, but the statement easily follows from their analysis). 

More recently, in joint work with Van Horn-Morris, the author was able to show that transverse knots classify contact structures in the following sense. Given a contact structure $\xi$ on a manifold $M$ and a knot $K$ in $M$, we let $\mathcal{T}^\xi_n(K)$ be the isotopy classes of transverse knots in the knot type $K$ with self-linking number $n$. 
\begin{theorem}[Etnyre and Van Horn-Morris 2011, \cite{EtnyreVanHorn-Morris11}]
Two contact structure $\xi_1$ and $\xi_2$ on $M$ are isotopic if and only if $|\mathcal{T}^{\xi_1}_n(K)|=|\mathcal{T}^{\xi_2}_n(K)|$ for all $n\in \Z$ and all fibered knots $K$ with pseudo-Anosov monodromy. 
\end{theorem}
We can paraphrase this theorem by saying the set of fibered knots with pseudo-Anosov monodromy \dfn{transversely characterizes} contact structures. More generally, if we can replace fibered knots with pseudo-Anosov monodromy in the above theorem with another class of knots $\mathcal{C}$, then we say $\mathcal{C}$ transversely classifies contact structures. We can analogously say that $\mathcal{C}$ \dfn{transversely characterizes tight contact structures} if the above theorem holds, where the $\xi_i$ are tight contact structures and $K$ is taken from $\mathcal{C}$. 

One can make a similar statement that Legendrian knots in a contact structure determine its isotopy class. While theoretically interesting, it is unlikely that the above theorem can be used to classify contact structures on a given manifold, but it does tell us that we can understand tight contact structures by understanding the transverse or Legendrian knot theory they support. Moreover, in known examples, a finite number of smooth knot types suffices to determine tight contact structures, and this fact has been used to classify tight contact structures. So it would be interesting to know if, given a manifold $M$, one could always find a finite collection of knots $\mathcal{C}$ that  %transversely (or Legendrianly) 
classify tight contact structures. %For example, three knots suffice for $T^3$.%For example, if $\mathcal{C}$ consists of any three linear knots in $T^3$, then $\mathcal{C}$ Legendrianly classifies tight contact structures on $T^3$. 

%%%%%%%%%%%%%%%%%%%%%%%%%%%%%%%%%%%%
\section{The structure of $1$-dimensional submanifolds in contact $3$-manifolds}
%%%%%%%%%%%%%%%%%%%%%%%%%%%%%%%%%%%%
While the general structure of transverse and Legendrian knots is still a mystery, we know quite a bit about their structure and have numerous results classifying Legendrian knots in a given knot type. It is known that the classification of Legendrian knots in a knot type determines the transverse classification, so we will focus mainly on Legendrian knots in this section. %While many of the results below hold in general tight contact manifolds, for the statements below, we will assume the knots are in the standard tight contact structure $\xi_{std}$ on $S^3$ (or $\R^3$), unless otherwise indicated. 

Given a knot type $K$ in a manifold $M$ and a contact structure $\xi$ on $M$, we let $\mathcal{L}(K)$ denote the set of isotopy classes of Legendrian knots smoothly isotopic to $K$. We have an obvious map
\[
\Phi\co \mathcal{L}(K)\to \Z\oplus \Z\co L\mapsto (\rot(L),\tb(L)).
\]
Classifying Legendrian knots in $\mathcal{L}(K)$ is equivalent to determining the image of $\Phi$, which is called the \dfn{mountain range of $K$}, and then the multiplicity of $\Phi$ for any point in the image. (Sometimes the mountain range of $K$ includes the multiplicity information and information about how the Legendrian knots are related by stabilization.) Determining the image of $\Phi$ is sometimes called the \dfn{geography problem} while determining the multiplicity of each point is sometimes called the \dfn{botany problem}. We will discuss what is known about these two problems below but begin by giving some classification results to give a sense as to what to expect.

%%%%%%%%%%%%%%%%%%%%%%%%%%%%%%%%%%%%
\subsection{Classification results}
%%%%%%%%%%%%%%%%%%%%%%%%%%%%%%%%%%%%
Here, we review what is known about the classification of Legendrian knots in various contact manifolds. 

%%%%%%%%%%%%%%%%%%%%%%%%%%%%%%%%%%%%
\subsubsection{Knots in $(S^3,\xi_{std})$}
\begin{figure}[htb]{
\begin{overpic}%[grid,tics=10] 
{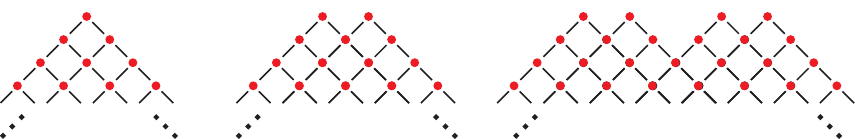}
\put(46, 61){\tiny$(0,t)$}
\put(105, 61){\tiny$(-1,-m-5)$}
\put(181, 61){\tiny$(1,-m-5)$}
%\put(360.5, 58){\color{red} \tiny$\mathbf a$}%272+88.5
%\put(349, 45){\color{red} \tiny$\mathbf b$}
%\put(371.5, 45){\color{red} \tiny$\mathbf b$}
%%
%\put(338, 34){\color{red} \tiny$\mathbf b$}
%\put(382, 34){\color{red} \tiny$\mathbf b$}
%%
%\put(326.5, 23){\color{red} \tiny$\mathbf b$}
%\put(393.5, 23){\color{red} \tiny$\mathbf b$}
\end{overpic}}
 \caption{On the left is the mountain range for the unknot (when $t=-1$), the figure eight knot (when $t=-3$), and the positive twist knots (when $t=-m-1$ for an even number $m$ of twists). Positive twist knots with an odd number $m$ of twists are shown in the middle. The right diagram is an example of the mountain range for a negative $(p,q)$-torus knots; the number of peaks and depth of the valleys depend on $p$ and $q$, and the peaks are at a height of $pq$. The diagonal lines indicate stabilizations. %On the right are the mountain ranges for the negative twist knots. When there are an odd number $m$ of twists $t=-3, a=-(m+1)/2,$ and $b=1$, while when there is an even number of twists $t=1, a=\lceil m^2/8\rceil,$ and $b=-\lceil m/2\rceil$. 
 }
  \label{fig:basicex}
\end{figure}

%%%%%%%%%%%%%%%%%%%%%%%%%%%%%%%%%%%%
We begin by recalling that a knot type $K$ is called \dfn{Legendrian simple} if the Legendrian isotopy class of a Legendrian knot in that knot type is determined by its rotation number and Thurston-Bennequin invariant, that is, if $\Phi$ is injective. In this case, we know that the classification of Legendrian representatives of $K$ is determined by its mountain range. Here are some examples of Legendrian simple knot types: the unknot \cite{EliashbergFraser98, EliashbergFraser09}, torus knots \cite{EtnyreHonda01b}, the figure eight knot \cite{EtnyreHonda01b}, cables of negative torus knots \cite{EtnyreHonda05}, and positive twist knots \cite{EtnyreNgVertesi13}. See Figure~\ref{fig:basicex}. 

The first smooth knot that was discovered not to be Legendrian simple was the twist knot $5_2$. In \cite{Chekanov02}, Chekanov showed that there were two Legendrian representatives of this knot that were not Legendrian isotopic, even though they have the same classical invariants. To distinguish these Legendrian knots, Chekanov defined a differential graded algebra (DGA) that was generated by the double points of the projection of the knot to the $xy$-plane and whose differential counted certain immersed polygons in the $xy$-plane with boundary on the projection. At the same time, Eliashberg, Givental, and Hofer \cite{EliashbergGiventalHofer00} were developing symplectic field theory, which for Legendrian knots in the standard contact $\R^3$ amounted to the above differential graded algebra, \cite{Eliashberg1998}. This DGA is now known as the \dfn{Chekanov-Eliashberg DGA}. 

Later, Ozsv\'ath, Szab\'o, and Thurston \cite{OzsvathSzaboThurston2008} and Lisca, Ozsv\'ath, Stipsicz, and Szab\'o \cite{LiscaOzsvathStipsiczSzabo2009} defined an invariant of Legendrian knots in knot Heegaard Floer homology. This invariant had the advantage that it could be non-zero for stabilized Legendrian knots, while the Chekanov-Eliashberg DGA is trivial for such knots. With this new invariant, one could show, among other things, that there were Legendrian realizations of twist knots that were stabilized but still distinct. 

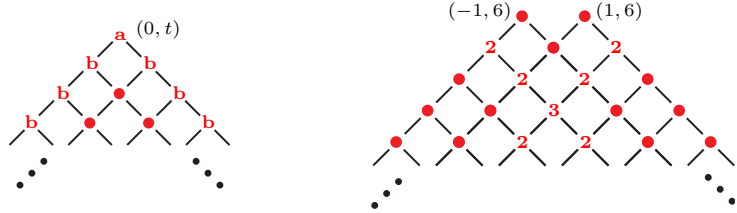
\begin{figure}[htb]{
\begin{overpic}%[grid,tics=10] 
{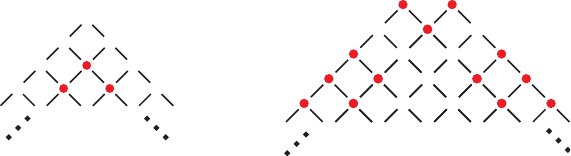}
\put(48, 66){\tiny$(0,t)$}
\put(40, 64){\color{red} \tiny$\mathbf a$}%360.5
\put(29, 53){\color{red} \tiny$\mathbf b$}
\put(51, 53){\color{red} \tiny$\mathbf b$}
\put(18, 41.5){\color{red} \tiny$\mathbf b$}
\put(62, 41.5){\color{red} \tiny$\mathbf b$}
\put(6, 30.5){\color{red} \tiny$\mathbf b$}
\put(73, 30.5){\color{red} \tiny$\mathbf b$}
\put(166, 73){\tiny$(-1,6)$}
\put(222, 73){\tiny$(1,6)$}
\put(180, 59){\color{red} \tiny$\mathbf 2$}
\put(227.5, 59){\color{red} \tiny$\mathbf 2$}
\put(192, 47){\color{red} \tiny$\mathbf 2$}
\put(215.5, 47){\color{red} \tiny$\mathbf 2$}
\put(192, 23){\color{red} \tiny$\mathbf 2$}
\put(216, 23){\color{red} \tiny$\mathbf 2$}
\put(204, 35){\color{red} \tiny$\mathbf 3$}
\end{overpic}}
 \caption{On the left is the mountain range for the negative twist knots. When there are $m$ twists and $m$ is odd, $t=-3, a=-(m+1)/2,$ and $b=1$, while when $m$ is even, $t=1, a=\lceil m^2/8\rceil,$ and $b=-\lceil m/2\rceil$. On the right is the mountain range for the $(3,2)$-cable of the right-handed trefoil.}
  \label{fig:nonsimp}
\end{figure}

Using the above two invariants and convex surface theory, Ng, V\'ertesi, and the author were able to completely classify Legendrian twist knots \cite{EtnyreNgVertesi13}. We noted above that the positive ones were Legendrian simple, but the negative twist knots are not. See Figure~\ref{fig:nonsimp}.% for their classification. 

We note that there are two other ways that non-Legendrian simple knot types have been found. One is via braid theory. In \cite{BirmanMenasco06II}, Birman and Menasco showed how to use a ``Markov theorem without stabilization" to prove that there are two non-transversely isotopic, transverse knots with the same self-linking number (which implies those knot types are not Legendrian simple). One may also use convex surface theory to distinguish and classify Legendrian knots. The first non-Legendrian simple knots found via this technique were due to Honda and the author \cite{EtnyreHonda05}, where they classified Legendrian knots in the knot type of the $(3,2)$-cable of the right-handed trefoil. See Figure~\ref{fig:nonsimp}. 
In \cite{EtnyreLafountainTosun12}, LaFountain, Tosun, and the author classified all cables of torus knots and showed that they had similar, but much more complicated behavior. For example, they were able to give the first examples of two Legendrian knots that had the same classical invariants, but were distinct after an arbitrarily large (but finite) number of positive and negative stabilizations. 

The above give specific classification results, but we also have several ``structure results". A structure result is a result that says how the classification of Legendrian knots is affected by some natural topological operations. The first such result was about connected sums of Legendrian knots. Honda and the author showed that the classification of Legendrian knots in the knot type of the connected sum $K_1\# K_2$ of knots $K_1$ and $K_2$ was determined by the classification in Legendrian knots realizing $K_1$ and $K_2$, \cite{EtnyreHonda03}. Specifically, there is a surjective map
\[
\Psi\co \mathcal{L}(K_1)\times \mathcal{L}(K_2)\to \mathcal{L}(K_1\# K_2)\co (L_1, L_2)\mapsto L_1\# L_2
\]
and the non-injectivity occurs only by sifting stabilizations from one of the knots to the other (or symmetries of $K_1\# K_2$). A surprising feature of this result is that the connected sum of Legendrian simple knots can produce non-Legendrian simple knots! For example, the connected sum of most negative torus knots results in non-Legendrian simple knot types. In fact, these were the first knot types for which we had a classification of Legendrian knots in a non-simple knot type.

In \cite{EtnyreHonda05, EtnyreLafountainTosun12} some progress was made on Legendrian representatives of cables of torus knots, but that heavily relied on the fact that torus knots were relatively simple. More recently, in \cite{ChakrabortyEtnyreMin2024Cabling}, Chakraborty, Min, and the author showed how to determine the classification of Legendrian knots in ``sufficiently large" cables of any knot in terms of the Legendrian classification of the knot; which, in particular, shows that a knot type is Legendrian simple if and only if any (or all) of its sufficiently large cables is (the forward implication was previously proven by Tosun, \cite{Tosun13}). They also gave a process to try to understand ``sufficiently negative" cables, and one can observe that such cables can be non-Legendrian simple even if the underlying knot is simple. Further progress was made understanding sufficiently negative cables in \cite{ChatterjeeEtnyreMinRodewald25}.

With connected sums and cables fairly well understood, the next interesting topological operation on a knot is taking satellites (cabling is a special case of this).  In \cite{EtnyreVertesi18}, V\'ertesi and the author gave a program to study satellite knot types. To understand Legendrian representatives of a satellite, one needs to understand, at a minimum, the Legendrian knots in the underlying knot type (frequently called the companion knot) as well as Legendrian representatives of the pattern knot (this is a knot in a solid torus). Even with this information, to obtain a complete classification of Legendrian representatives, one also needs to know that the companion knot has a technical condition called uniform thickness. 

%%%%%%%%%%%%%%%%%%%%%%%%%%%%%%%%%%%%
\subsubsection{Legendrian links in $(S^3,\xi_{std})$}
%%%%%%%%%%%%%%%%%%%%%%%%%%%%%%%%%%%%
Legendrian representatives of links have been studied much less than for knots. Here, apart from classification results, one can ask for structural results. The two most natural questions along these lines are a geography question and a symmetry question. The geography question asks if given Legendrian realizations of the components of a link can they be realized by a Legendrian realization of the link. The symmetry question asks what topological symmetries of a knot can be realized by Legendrian isotopies. 

One of the first geography results was due to Mohnke. In \cite{Mohnke01}, he used a refinement of the Bennequin inequality to show that some links whose components were unknots could not be realized by a Legendrian link all of whose components were maximal Thurston-Bennequin invariant Legendrian unknots. 

The first symmetry result was due to Traynor. In \cite{Traynor97}, she used generating function techniques to show that there is a Legendrian two-component link in the jet space of $S^1$ whose components could not be interchanged via a Legendrian isotopy even though they could be interchanged by a smooth isotopy. The first symmetry result in $S^3$ was due to Mishachev. He used the Chekanov-Eliashberg DGA to show that one could only cyclically permute the components of the ``$n$-copy" of the maximal Thurston-Bennequin unknot via a Legendrian isotopy, \cite{Mishachev03}. The next step in understanding what is possible for symmetries of a Legendrian knot can be found in \cite{Ng03}, where Ng used the Chekanov-Eliashberg DGA to show that no components of the $n$-copy of the maximal Thurston-Bennequin invariant Legendrian figure-eight knot could be permuted via a Legendrian isotopy. 

Another surprising example of Legendrian links was provided by Jordan and Traynor in \cite{JordanTraynor06}. There, they showed that there were Legendrian links that were smoothly isotopic, and each component of one link was Legendrian isotopic to a component of the other link, but the links were not Legendrian isotopic. 

The first classification result for Legendrian links was due to Ding and Geiges. In \cite{DingGeiges07}, they classified what they called ``cable links". These are links consisting of an unknot and a $(p,q)$ torus knot sitting on the boundary of a neighborhood of the unknot. When $p>1$, these links are Legendrian simple and there are no topological symmetries. In \cite{DingGeiges10}, they were able to show that when $p=1$ and we consider the link in the $1$-jet space of $S^1$, then there are no Legendrian isotopies that interchange the components, thus giving the first classification result where smooth and Legendrian symmetries were different (as noted above, it was previously known that these types of symmetries were different, but only examples were know, and not a complete classification). In \cite{GeigesOnaran20a}, Geiges and Onaran classified Legendrian realizations of the Hopf link. 

In \cite{DaltonEtnyreTraynor2024LegLink}, Dalton, Traynor, and the author classified Legendrian torus links and also Legendrian representatives of cable links of Legendrian simple and uniformly thick knot types. In this work, they showed that there are serious restrictions on which Legendrian knots could be realized as the components of a Legendrian link, that there are non-destabilizable representatives of Legendrian links that do not have maximal Thurston-Bennequin number, and that some cable links have no Legendrian symmetries, others have only cyclic symmetries, while still others have the same Legendrian symmetries as smooth symmetries. Then, in \cite{ChatterjeeEtnyreMinRodewald25}, Chatterjee, Min, Rodewald, and the author extended these results to cable links for non-simple knot types.  Doing so shows that there can be two Legendrian links that are smoothly isotopic, have highly stabilized components that are Legendrian isotopic, but are not Legendrian isotopic as links. We were also able to show that for a uniform thick knot type, any two components of a Legendrian realization of a cable link with the same classical invariants must be isotopic. This severely limits the Legendrian knots that can be realized as components of a Legendrian cable link. 

%%%%%%%%%%%%%%%%%%%%%%%%%%%%%%%%%%%%
\subsubsection{Legendrian links in other tight contact manifolds}
%%%%%%%%%%%%%%%%%%%%%%%%%%%%%%%%%%%%
There has not been a great deal of study of Legendrian knots in manifolds other than $S^3$, but there have been a few results. We first note that in most cases, classification results are ``coarse classification results" in the sense that Legendrian knots are only determined up to contactomorphism (smoothly isotopic to the identity) instead of up to isotopy. %as in the case for the tight contact structure on $S^3$. %This is due to the fact that we do not have Eliashberg's result \cite{Eliashberg92a} about contactomorphism of the tight contact structure on $S^3$ for a general contact manifold. 

The first classification result was due to Onaran. In \cite{Onaran18} she showed that ``positive" torus knots in lens spaces with their universally tight contact structure were Legendrian simple and the possible classical invariants were determined. Here, a torus knot means a knot sitting on a Heegaard torus for the lens space and we refer to the paper for the definition of ``positive". In \cite{ChenDingLi15}, Chen, Ding, and Li showed that torus knots in the unique tight contact structure on $S^1\times S^2$ were Legendrian simple. 
%(We note that these knots are not null-homotopic, so one must define the Thurston-Bennequin invariant in a different way.) 
In \cite{BakerEtnyre12} Baker and the author showed that rational unknots, which are cores of Heegaard tori, in tight contact structures on lens spaces, were Legendrian simple and determined their possible Thurston-Bennequin invariants, and for $L(p,1)$ with $p$ odd, they also gave their possible rotation numbers, completing the classification. They also indicated how this worked for a general $L(p,q)$. In \cite{GeigesOnaran15}, Geiges and Onaran explicitly worked out the case for $L(p,1)$ when $p$ is even as well as the case for Legendrian knots in $L(5,2)$. 

%%%%%%%%%%%%%%%%%%%%%%%%%%%%%%%%%%%%
\subsection{The structure of Legendrian knots}
%%%%%%%%%%%%%%%%%%%%%%%%%%%%%%%%%%%%
In this section, we discuss the general structure of Legendrian knots in a given knot type. As mentioned at the beginning of this section, we can break this discussion into two parts: the geography of Legendrian knots, that is, the image of $\Phi$ for a given knot type, and the botany of Legendrian knots, that is, the preimage of a point under $\Phi$. 

We begin with the botany problem. Here are the general things we know. %about the botany of Legendrian knots. %If you fix knot type $K$ then we can say the following about the preimage of points under $\Psi$:
\begin{enumerate}
\item For any $(r,t)\in\Z^2$, $\Phi^{-1}(r,t)$ is finite, \cite{ColinGirouxHonda09}.
\item However, $\Phi^{-1}(r,t)$ can be arbitrarily large, \cite{EtnyreHonda03, EtnyreLafountainTosun12, EtnyreNgVertesi13}.
\item Any two elements in $\Phi^{-1}(r,t)$ are isotopic after sufficiently many positive and negative stabilizations, %many positive and negative stabilziations, 
\cite{FuchsTabachnikov97}. 
\item There are elements in $\Phi^{-1}(r,t)$ that  
\begin{enumerate}
\item can be positively (or negatively) stabilized infinitely often and remain distinct, \cite{EtnyreHonda03, EtnyreHonda05, EtnyreLafountainTosun12, EtnyreNgVertesi13}, and
\item require arbitrarily many stabilizations of both signs before they become isotopic, \cite{EtnyreHonda03, EtnyreLafountainTosun12}. 
\end{enumerate}
\end{enumerate}

Turning to the geography of Legendrian knots, we know the following:
\begin{enumerate}
\item In tight $S^3$, the image of $\Phi$ is easily seen to be symmetric, but in other manifolds it might not be, \cite{BakerEtnyre12, GeigesOnaran15}. 
\item The sum of the coordinates of any point in the image of $\Phi$ must be odd. 
\item If $(r,t)$ is in the image of $\Phi$, then using stabilizations we see that the cone $C((r,t))=\{(r,t-n): n\geq 0, j=-n,-n+2,\ldots, n-2, n \}$ is also in the image. 
\item The image of $\Phi$ is bounded above, \cite{Bennequin83, FuchsTabachnikov97}.
\item The upper bounds in  \cite{Bennequin83, FuchsTabachnikov97} can be arbitrarily far from being sharp, \cite{EtnyreHonda01b}.
\item The mountain range can have arbitrarily many peaks and deep valleys, \cite{EtnyreHonda01b}.
\item There can be non-maximal Thurston-Bennequin elements that do not destabilize, that is, the multiplicity of the points in the mountain range can increase as one moves to smaller values of $\tb$, \cite{EtnyreHonda05, EtnyreLafountainTosun12}. 
\item There can be non-maximal peaks, \cite{ChongchitmateNgPre}. 
\end{enumerate}
We note that for the last item, we do not have a classification of Legendrian knots in a knot type where we have non-maximal peaks, so we do not have a good understanding of why they occur, but only that they do. 

From above, we have a good sense of the general structure of the classification of Legendrian knots, but there are some notable things we do not know. 
\begin{question}
Given a smooth knot type $K$, is there an integer $n$ such that a Legendrian knot $L\in\mathcal{L}(K)$ with $\tb(L)<n$ destabilizes?
\end{question}
We note that if the answer to this question is ``Yes", then $\mathcal{L}(K)$ is finitely generated in the following sense. Given the upper bound on the image of $\Phi$, we see that there are finitely many points with second coordinate larger than $n$ (this is the $n$ in the question). Then the finiteness result for $\Psi^{-1}(r,t)$ tells us that there are finitely many Legendrian knots with $\tb\geq n$. Thus, all elements in $\mathcal{L}(K)$ are stabilizations of a finite number of non-detabilizable elements. So any classification result amounts to identifying these finitely many elements and understanding when they become the same under stabilizations. 

\begin{question}
Are there points in the interior of the image of $\Phi$ that do not destabilize?
\end{question}
This question does not have the impact of the previous question, but it would help us understand the overall behavior of Legendrian representatives of a given knot type. 

%%%%%%%%%%%%%%%%%%%%%%%%%%%%%%%%%%%%
\subsection{Knots in overtwisted contact structures}
%%%%%%%%%%%%%%%%%%%%%%%%%%%%%%%%%%%%
While Legendrian knots in overtwisted contact structures have recently become an important area of research, there has historically not been a great deal of study of them. This is likely due to the fact that Eliashberg has classified overtwisted contact structures \cite{Eliashberg89}, and they seemed relatively simple compared to tight contact structures. Moreover, most of the applications of contact geometry to low-dimensional topology are connected to tight contact structures. However, 
in recent years, it has become clear that there is not only a surprisingly rich theory of Legendrian knots in overtwisted contact manifolds, but they can be instrumental in understanding the classification of tight contact structures on $3$-manifolds. We will discuss the first point here and the second in the next subsection. 

The classification of Legendrian knots in overtwisted contact structures suffers from not having the use of the Chekanov-Eliashberg DGA, augmentations, or other holomorphic tools used in the tight case. However, we still have the invaluable tools from Heegaard Floer theory \cite{LiscaOzsvathStipsiczSzabo2009}. We also still have the use of convex surface theory, the main tool used in almost all classification results, but due to the nature of overtwisted contact structures, many of the arguments are more difficult and require new techniques. 

We first note that Legendrian knots $L$ in an overtwisted contact structure $\xi$ come in one of two types. We say $L$ is \dfn{loose} if the contact structure on its complement is overtwisted, and otherwise we call it \dfn{non-loose} (some authors call such knots \dfn{exceptional}). It is not hard to use Eliashberg's result mentioned above to show that loose knots are determined by their classical invariants and that any pair of integers $(r,t)$ for which $r+t$ is odd, can be realized as the rotation number and Thurston-Bennequin invariant of a loose Legendrian knot, \cite{Etnyre13}. Once again, we note that ``classification" here means the knots are determined up to contactomorphism smoothly isotopic to the identity.%, this is usually called a \dfn{coarse classification}. 

The first general existence result for non-loose Legendrian representatives of a knot appeared in \cite{Etnyre13}, but there was an error in the author's paper that was corrected in the work of Chatterjee, Mukherjee, Min, and the author. In \cite{ChatterjeeEtnyreMinMukherjee25nonloose}, they showed that a knot $K$ in an irreducible $3$-manifold $M$ admitted a non-loose Legendrian representative if and only if either $K$ is not contained in a $3$-ball, or $M$ admits a tight contact structure. So we know that non-loose Legendrian knots exist in abundance, but, as we will see below, it seems that they only exist in very few overtwisted contact structures on a given manifold. 

Most of the classification results we discuss will be in $S^3$ with some overtwisted contact structure on it. We note that there are an integer's worth of such contact structures \cite{Eliashberg89, Eliashberg92a}, so we denote them by $\xi_n$. Here $n$ is the so-called $d_3$-invariant of the plane field.

\begin{figure}[htb]{
\begin{overpic}%[grid,tics=10] 
{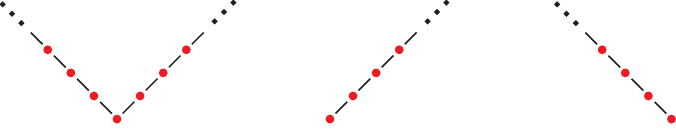}
%\put(46, 61){\tiny$(0,t)$}
\end{overpic}}
 \caption{On the left, we see the mountain range for the non-loose unknots in $S^3$ where the bottom vertex is at $(0,1)$. All three figures can be the mountain range for a non-loose rational unknot in some lens space. The coordinates of the bottom vertex and whether one has a ``V" as ``forward slash" or ``back slash" depend on the lens space and contact structure. }
  \label{fig:nl}
\end{figure}

The first classification result for non-loose knots was due to Eliashberg and Fraser. In \cite{EliashbergFraser09}, they showed that non-loose Legendrian unknots exist only in $\xi_1$ and there they are coarsely classified by their classical invariants. Figure~\ref{fig:nl} indicates the possible values for these invariants. 

Not surprisingly, many of the structural features of the mountain range for Legendrian knots in tight contact manifolds do not hold here. Most obviously, we do not have a Bennequin inequality, but we do have the following inequality due to \'{S}wi\k{a}tkowski, \cite{Swiatkowski92}, for the invariants of a non-loose Legendrian knot
\[
-|\tb(L)|+|\rot(L)|\leq -\chi(\Sigma)
\]
where $\Sigma$ is a surface with boundary $L$. 

The analogs for unknots in lens spaces are called rational unknots and are the cores of Heegaard tori in the lens spaces. Geiges and Onaran \cite{GeigesOnaran15} classified non-loose rational unknots in $L(p,1)$ and also in $L(5,2)$. In \cite{ChatterjeeEtnyreMinMukherjee25nonloose} Chatterjee, Min, Mukherjee, and the author completed the classification in all lens spaces and showed that the mountain range of non-loose knots in a fixed overtwisted contact structure was of the form shown in Figure~\ref{fig:nl}.

When considering non-loose knots, one can ask if they have Giroux torsion in their complement (this is a particular contact structure on a thickened torus). The next set of results only considered those that do not contain Giroux torsion in their complement; these are sometimes called \dfn{strongly exceptional}. In \cite{GeigesOnaran20}, Geiges and Onaran classified non-loose Legendrian representatives of the right and left-handed trefoils with specific values for the classical invariants. They were also able to classify $(p,-(pn-1))$-torus knots with certain classical invariants. Later, Matkovi\v{c} \cite{Matkovic22}, was able to classify non-loose Legendrian negative $(p,q)$-torus knots with no Giroux torsion in their complement when their Thurston-Bennequin invariant was less than $pq$. 

\begin{figure}[htb]{
\begin{overpic}%[grid,tics=10] 
{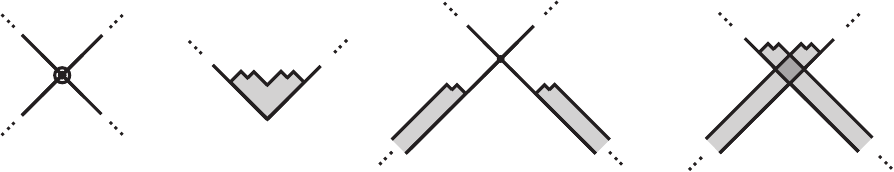}
%\put(46, 61){\tiny$(0,t)$}
\end{overpic}}
 \caption{Possible mountain ranges for non-loose Legendrian torus knots in $S^3$. }
  \label{fig:nlt}
\end{figure}

In \cite{EtnyreMinMukherjee22Pre}, Min, Mukherjee, and the author were able to give a complete classification of non-loose Legendrian torus knots with or without Giroux torsion in their complement. The classification is expressed in terms of a simple algorithm, and several families of examples were written out in closed form. From this work, we can observe a few interesting facts. First, there are only a finite number of overtwisted contact structures on $S^3$ that contain non-loose Legendrian $(p,q)$-torus knots, and the number of such contact structures can be determined by the continued fraction of $p/q$. Secondly, the mountain range for non-loose Legendrian knots in any fixed overtwisted contact structure is of the form shown in Figure~\ref{fig:nlt}. 

Turning to links, Geiges and Onaran \cite{GeigesOnaran20a} classified non-loose Hopf links in $S^3$. Then in \cite{ChatterjeeGeigesOnaran24pre}, they, together with Chatterjee, classified non-loose Hopf links in the overtwisted contact structures on $L(p,1)$, while in forthcoming work, Chatterjee will complete the classification in all lens spaces. 

We end this section by discussing a particularly interesting example where non-loose knots have been classified up to Legendrian isotopy. In the groundbreaking paper \cite{Vogel2018}, Vogel determined the contact mapping class group of overtwisted contact structures on $S^3$, thus allowing him to classify non-loose Legendrian unknots up to isotopy. Specifically, he showed that each of the pairs $(r,t)$ realized by the Eliashberg-Fraser result as the rotation number and Thurston-Bennequin invariant of a non-loose Legendrian is realized by two distinct Legendrian knots, up to Legendrian isotopy. He also showed that loose Legendrian unknots in $\xi_1$ are determined by their classical invariants if and only if their Thurston-Bennequin invariant is less than $0$.

%%%%%%%%%%%%%%%%%%%%%%%%%%%%%%%%%%%%
\subsection{Constructing and classifying contact structures via Legendrian knots}
%%%%%%%%%%%%%%%%%%%%%%%%%%%%%%%%%%%%
Legendrian knots have always been one of the main ways to construct tight contact structures via surgery in symplectically fillable contact structures, but in recent years, it has become clear that they can also be instrumental in classifying tight contact structures as well. Consider a manifold $Y$ obtained by Dehn surgery on a knot $K'$ from another manifold $Y'$ (frequently $Y'$ will be $S^3$ or some fairly simple $3$-manifold). Then inside the surgery torus in $Y$ there is a dual knot $K$, so the appropriate surgery on $K$ will produce $Y'$. Given any tight contact structure $\xi$ on $Y$, we can take a Legendrian knot $L$ in the knot type of $K$. We can take a neighborhood of $L$ and then expand it until it is the ``largest contact neighborhood possible" (this has a precise definition, see \cite{EtnyreMinMukherjee22Pre}). Now the complement of a neighborhood of this $L$ will be a contact structure on the complement of a neighborhood of $K'$ in $Y'$. In many situations, we can assume that the contact structure on the complement of $K'$ defines a Legendrian knot in $Y'$ with some contact structure $\xi'$, which might be overtwisted. In our setup, this will either be a Legendrian knot in a tight contact structure on $Y'$ or a non-loose Legendrian knot in an overtwisted contact structure on $Y'$.  (It turns out there is another case, that of ``non-thickenable" tori, see \cite{EtnyreLafountainTosun12}, but these are naturally found when studying Legendrian knots in contact structures, so we will elide these to simplify the discussion.) Thus, any tight contact structure on $Y$ is obtained from some contact surgery on a Legendrian representative of $K'$ in some contact structure on $Y'$. We can then use Heegaard Floer invariants to show that all the contact structures obtained by such contact surgeries are tight on $Y$ and all are different. (In the elided case, we must use more geometric arguments.)

In the upcoming paper \cite{EtnyreMinTosunVarvarezosPre} by Min, Tosun, Varvarezos, and the author, the above procedure is being carried out to classify tight contact structures on manifolds obtained by any surgery on a torus knot in $S^3$. Currently, the only classification of contact structures on manifolds obtained by surgery on a knot is for the unknot. So this result now gives such a result on an infinite family of knots. Moreover, many of the special features of contact structures have been observed on manifolds obtained by surgery on torus knots (for example, the only manifolds known not to admit tight contact structures come from such surgeries \cite{EtnyreHonda01a, LiscaStipsicz07}, and the first, and most, tight but not symplectically fillable contact structures also appear from such surgeries \cite{EtnyreHonda02b, LiscaStipsicz04}). 

Min, Xu, and the author are currently classifying non-loose Legendrian torus knots in lens spaces, and using this, Min, Tosun, Varvarezos, and the author plan to complete the classification of tight contact structures on small Seifert fibered spaces. Such a classification has been a driving goal for around 25 years and has led to the discovery of many new phenomena in contact geometry. Having this complete classification will create a test bed for conjectures and the further study of tight contact structures in dimension 3. 

In addition, Min, Le, and the author are working on the classification of non-loose Legendrian representations of the figure eight knot, which in turn is expected to lead to a classification of tight contact structures on manifolds obtained from surgery on the figure eight knot. This would be the first such classification for a hyperbolic knot, though Conway and Min have already made substantial progress on this, \cite{ConwayMin20}. At this point, it appears that we will see new phenomena for non-loose Legendrian knots as well as new phenomena for tight contact structures on hyperbolic manifolds. 

%%%%%%%%%%%%%%%%%%%%%%%%%%%%%%%%%%%%
\section{Submanifolds in higher dimensions}
%%%%%%%%%%%%%%%%%%%%%%%%%%%%%%%%%%%%
The theory of contact structures in high dimensions is not as well developed as it is in dimension $3$. But there have been some notable results in the last 10 to 15 years. Some high points are Casals, Pancholi, and Presas's proof that all almost contact $5$-manifolds admit a contact structure, \cite{CasalsPancholiPresas2015}, Borman, Eliashberg, and Murphy's definition of overtwisted contact structures and proof that all almost contact manifolds admit a contact structure \cite{BormanEliashbergMurphy15}, and Bowden, Gironella, Moreno, and Zhou's proof that tight but not symplectically fillable contact structures are plentiful \cite{BowdenGironellaMorenoZhou24Pre}.

Below, we will discuss the role of submanifolds in high-dimensional contact geometry. We begin by seeing how to construct contact structures via surgery on special submanifolds. We then turn to distinguishing special submanifolds, starting with Legendrian submanifolds in Section~\ref{DGA} and then contact submanifolds in Section~\ref{noniso}. 

%%%%%%%%%%%%%%%%%%%%%%%%%%%%%%%%%%%%
\subsection{Constructing contact structures}
%%%%%%%%%%%%%%%%%%%%%%%%%%%%%%%%%%%%
Given a contact manifold $(M,\xi=\ker \alpha)$ of dimension $2n+1$, a submanifold $S$ is called \dfn{isotropic} if it is tangent to $\xi$; this will imply that the dimension of $S$ is at most $n$ and when it has dimension $n$ it will be called a \dfn{Legendrian submanifold}. The submanifold is called \dfn{co-isotropic} if the $d\alpha$ orthogonal complement of $T_xS\cap \xi_x$, which we denote by $S_\xi$, in $\xi$ is contained in $S_\xi$. It has long been known how to perform a contact surgery on an isotropic sphere (with a suitable trivialization of its normal bundle) \cite{Weinstein91}, and in \cite{ConwayEtnyre20caps}, Conway and the author discussed a contact surgery on a co-istoropic sphere, and proved the following result.
\begin{theorem}[Conway and Etnyre 2020, \cite{ConwayEtnyre20caps}]
Any oriented $2n+1$-dimensional closed contact manifold $(M,\xi)$ can be obtained from the standard contact structure on $S^{2n+1}$ by contact surgeries on isotropic and co-isotropic spheres. 
\end{theorem}
We note that Lazarev also established this result for ``almost-Weinstein fillable contact manifolds" using a subset of the surgeries we used, \cite{Lazarev2020caps}. A nice corollary of the above theorem, and Lazarev's work, is that any contact manifold admits infinitely many distinct symplectic caps. A symplectic cap is a compact symplectic manifold with the contact manifold as its concave boundary. Another interesting corollary, first observed by Lazarev but proved in both their papers, is that given any finite set of contact manifolds, there is one fixed contact manifold to which they are all symplectically cobordant. 

%%%%%%%%%%%%%%%%%%%%%%%%%%%%%%%%%%%%
\subsection{Distinguishing Legendrian submanifolds}\label{DGA}
\begin{figure}[htb]{
\begin{overpic}%[grid,tics=10] 
{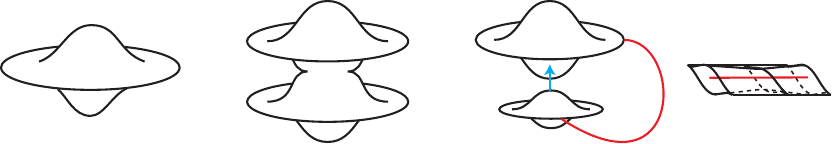}
%\put(46, 61){\tiny$(0,t)$}
\end{overpic}}
 \caption{The projections from $\R^5$ to $\R^3$ of three Legendrian spheres distinguished by the DGA. On the left, the tube about the red path is added along the curved path next to it, and the bottom part of the picture is pushed along the arrow to ``stabilize" the Legendrian sphere.}
  \label{fig:hdex}
\end{figure}

%%%%%%%%%%%%%%%%%%%%%%%%%%%%%%%%%%%%
Legendrian submanifolds are the most studied submanifolds in high-dimensional contact geometry. Though there were some nice constructions, there were few tools to distinguish them until 2005. Inspired by Chekanov and Eliashberg's work discussed above, and the general framework of symplectic field theory, Ekholm, Sullivan, and the author defined an analog of the Chekanov-Eliashberg DGA for Legendrian knots in the standard contact structure on $\R^{2n+1}$ in \cite{EkholmEtnyreSullivan05a, EkholmEtnyreSullivan05b, EkholmEtnyreSullivan05c} and were able to show that it is quite effective in distinguishing Legendrian submanifolds. See Figure~\ref{fig:hdex}. They later generalized this work to Legendrian submanifolds of the product of an exact symplectic manifold with $\R$ in \cite{EkholmEtnyreSullivan07}.

Since then, several authors have developed interesting structures for this DGA, relations with other ``homology theories", as well as nice applications. See, for example, \cite{BourgeoisEkholmEliashberg12, EkholmEtnyreSabloff09, Dimitroglou2016}. We also would like to point out that many of these results can also be approached from the perspective of generating family homology \cite{SabloffTraynor13} and micro-local sheaves \cite{CasalsZaslow2022, GuillermouKashiwaraSchapira2012, TreumannZaslow2018}.% and micro-local sheaf theory \cite{}. {\color{red} Find good references}

We end this section by discussing a surprising application of this work to the topology of smooth knots in $\R^3$. It is easy to see that the unit cotangent bundle of $\R^3$ is a contact manifold and actually contactomorphic to the $1$-jet space of $S^2$, $J^1(S^2)$, and hence the DGA invariant discussed above is well defined for Legendrian submanifolds. Given a smooth knot $K$ in $\R^3$, we can consider its unit co-normal bundle $L_K$, which is the set of unit co-vectors in the cotangent bundle that vanish on the knot. It turns out that $L_K$ is a Legendrian torus in $J^1(S^2)$ and an isotopy of $K$ produces a Legendrian isotopy of this torus. Thus, the DGA associated with $L_K$ is an invariant of the smooth knot type of $K$. So we can use contact geometry to study knot theory! The homology of the DGA is usually called the \dfn{knot contact homology} of $K$. 

Inspired by this connection, Ng \cite{Ng05, Ng05a} wrote down a combinatorial invariant of $K$ that was expected to be this DGA and showed that it was a powerful invariant of knots. For example, it could distinguish the unknot from other knots (always a basic test of the strength of a knot invariant) and was related to other interesting knot invariants defined in purely topological ways. Moreover, in \cite{GordonLidman2017KnotContactHomology}, Gordon and Lidman showed that the invariant can determine if a knot is cabled, composite, or a torus knot. In \cite{EkholmEtnyreNgSullivan13}, Ekholm, Ng, Sullivan, and the author were able to show that Ng's invariant was indeed knot contact homology. So, knot contact homology is in fact a powerful invariant of knots. 

More recently, Shende \cite{Shende2019CompleteInvariant} was able to show that $L_K$ was a complete invariant of $K$. That is, $K$ and $K'$ are smoothly isotopic if and only if $L_K$ and $L_{K'}$ are Legendrian isotopic. Shende used microlocal sheaf theory to prove this. However, in \cite{EkholmNgShende18}, Ekholm, Ng, and Shende were able to show that a slight enhancement of knot contact homology could also distinguish any two smooth knots. It is an interesting, ongoing problem to determine how to understand important topological features of a knot in terms of its knot contact homology. 

We end this section by noting that one can enhance the knot contact homology to study transverse knots in the standard contact structure on $\R^3$. Specifically, one can add a filtration to the knot contact homology DGA for $L_T$ when $T$ is a transverse knot. This was carried out by Ekholm, Ng, Sullivan, and the author in \cite{EkholmEtnyreNgSullivan13a}. While in \cite{Ng2011transversehomology}, Ng showed that this was a very effective invariant of transverse knots, though we do not know if it is a complete invariant. 

It is ongoing work to study this ``conormal construction" for submanifolds in other dimensions and codimensions. 

%%%%%%%%%%%%%%%%%%%%%%%%%%%%%%%%%%%%
\subsection{Contact submanifolds}\label{noniso}
%%%%%%%%%%%%%%%%%%%%%%%%%%%%%%%%%%%%
We now turn to ``transverse submanifolds" in contact manifolds. One could mean several things by this, but the most geometric version of this is to consider submanifolds $S$ in a contact manifold $(M,\xi)$ such that $\xi\cap TS$, which we denote $\xi_S$, is a contact structure on $S$. We will call such an $S$ a \dfn{contact submanifold}. Two basic questions about contact submanifolds are: (1) Which contact manifolds can be embedded in a higher-dimensional contact manifold? That is, which contact manifolds can be realized as a contact submanifold? and (2) Given two contact submanifolds, when are they isotopic through contact submanifolds? We will survey what is known about each of these questions in the following two subsections. 

%%%%%%%%%%%%%%%%%%%%%%%%%%%%%%%%%%%%
\subsubsection{Existence of contact submanifolds}
%%%%%%%%%%%%%%%%%%%%%%%%%%%%%%%%%%%%
The theory of contact embeddings of $S^1$ into a contact $3$-manifold is the theory of transverse knots that we have discussed above. In higher dimensions, contact embeddings satisfy an $h$-principle \cite{Gromov86} if the codimension of the embedding is larger than 2, and thus they exist in abundance. However, codimension $2$ contact embeddings are more difficult to understand (though we will see that ultimately, they do satisfy an existence $h$-principle). 

One of the first general results was due to Kasuya, who in \cite{Kasuya13}, observed that the first Chern class of a contact structure that has a codimension $2$ contact embedding into an oriented contact manifold with trivial second cohomology must vanish. This, for example, excluded many contact $3$-manifolds from embedding in the standard contact structure $\xi_{std}$ on $S^5$. However, in \cite{EtnyreFurukawa17}, Furukawa and the author developed the notion of a braided embedding to show that many contact $3$-manifold with trivial first Chern class did indeed embed in $(S^5,\xi_{std})$. They also showed that if an embedding of a $3$-manifold $M$ into $S^5$ could be isotoped to a braided embedding, then it could be isotoped to be a contact embedding for some contact structure on $M$. Then in \cite{EtnyreLekili18}, Lekili and the author used open book embeddings to show there was a simple contact $5$-manifold, namely a Stein fillable contact structure on the twisted $S^3$-bundle over $S^2$, into which all contact $3$-manifolds embed. 

Finally, in \cite{CasalsPancholiPresas2021Whitney}, Casals, Pancholi, and Presas developed a Legendrian version of the Whitney trick to show that there is a general existence $h$-principle for codimension $2$ contact embeddings. Specifically, if there is a ``formal embedding" of a contact manifold into a higher-dimensional contact manifold, then there is actually a contact embedding. Here ``formal embedding" just means that there is an embedding whose differential is isotopic, as a bundle map, to a map on the tangent space that maps one of the contact structures into the other. This answers the general existence question completely, but there were two results appearing at a similar time that, taken together, also manage to establish this. Namely, in \cite{HondaHuang2023PreConvexHS}, Honda and Huang developed the theory of convex hypersurfaces in high dimensions and used it to prove that an ``almost contact" embedded codimension $2$ manifold in a contact manifold could be approximated by a contact embedding of some contact structure on the manifold; while Pancholi and Pandit \cite{PancholiPandit2022IsoContact} showed that given a contact embedding of some contact structure on a manifold then any other contact structure in the same almost contact type can also be contact embedded.

%%%%%%%%%%%%%%%%%%%%%%%%%%%%%%%%%%%%
\subsubsection{Contact submanifolds}
%%%%%%%%%%%%%%%%%%%%%%%%%%%%%%%%%%%%
Gromov's $h$-principle mentioned above also tells us that two contact embeddings of codimension larger than $2$ that are ``formally isotopic" are actually isotopic through contact embeddings, so the main question remaining concerns codimension $2$ contact embeddings. The first progress on this occurred only recently. In \cite{CasalsEtnyre2020NonIsotopic}, Casals and the author gave the first examples of two codimension $2$ contact embeddings that were formally isotopic, but not isotopic through contact embeddings. To accomplish this, they introduced the notion of the contact push-off of a Legendrian submanifold. This process produces a codimension $2$ contact submanifold associated with a Legendrian knot. The examples were then constructed by taking the contact push-off of different Legendrian knots, which were formally Legendrian homotopic, and they distinguished them by considering the contact structure induced on the $2$-fold branched cover of the contact push-offs. An interesting open question is whether or not the contact push-off of a Legendrian submanifold determines the Legendrian submanifold. %This is even open in dimension $3$. 

The next step in the study of codimension $2$ embeddings was due to Zhou. In \cite{Zhou2023Nonisotopic}, it was shown that in dimensions $7$ and above, one could use the construction above to distinguish infinitely many different codimension $2$ contact embeddings. More recently, C\^ot\'e and Fauteux-Chapleau \cite{CoteFauteux-Chapleau2024InvariantOfEmbedding} showed how to build a filtered version of contact homology to define an invariant of codimension $2$ contact submanifolds and then showed that this invariant was indeed effective at distinguishing formally isotopic contact embeddings. 

In a recent preprint, Avdek \cite{Avdek2024Pre} defined a notion of stabilization for codimension $2$ contact submanifolds. From this, he could also give examples of non-contact isotopic embeddings that were formally isotopic. This is the first general modification of contact embeddings and promises to be an important direction for the exploration of contact embeddings. 

%%%%%%%%%%%%%%%%%%%%%%%%%%%%%%%%%%%%
\subsection{Conclusion}
%%%%%%%%%%%%%%%%%%%%%%%%%%%%%%%%%%%%
We hope this note has provided a good historical overview of the deep connections between the study of submanifolds of contact manifolds and the development of the rich and beautiful subject of contact geometry. While there has been, as noted above, much groundbreaking work in the field, there is still a great deal of work to be done and many important questions to be answered and many more to be discovered. 

\def\cprime{$'$}

%% references
%\bibliography{references}
%%\bibliographystyle{alpha}
%\bibliographystyle{plain}
\end{document}